\documentclass{amsart}
\usepackage{amssymb,latexsym,amsmath,amsfonts,hhline}
\usepackage{pdfsync}
\usepackage{xcolor} 

\makeatletter
\renewcommand{\subsection}{\@startsection{subsection}{2}{0mm}{-2mm}{-2mm}{\bf\normalsize}}

\def\sbsnt#1{\subsection{#1}}
\makeatother

\newtheorem{formula}{}[section]
\newtheorem{definition}[formula]{Definition}
\newtheorem{corollary}[formula]{Corollary}
\newtheorem{remark}[formula]{Remark}
\newtheorem{lemma}[formula]{Lemma}
\newtheorem{theorem}[formula]{Theorem}

\def\thrm{\begin{theorem}}
\def\thrml#1{\begin{theorem}\label{#1}}
\def\ethrm{\end{theorem}}
\def\rmrk{\begin{remark}}
\def\rmrkl#1{\begin{remark}\label{#1}}
\def\ermrk{\end{remark}}
\def\dfntn{\begin{definition}}
\def\dfntnl#1{\begin{definition}\label{#1}}
\def\edfntn{\end{definition}}
\def\nmrt{\begin{enumerate}}
\def\enmrt{\end{enumerate}}
\def\tm#1{\item[{\rm (#1)}]}
\def\qtn{\begin{equation}}
\def\qtnl#1{\begin{equation}\label{#1}}
\def\eqtn{\end{equation}}
\def\lmm{\begin{lemma}}
\def\lmml#1{\begin{lemma}\label{#1}}
\def\elmm{\end{lemma}}
\def\crllr{\begin{corollary}}
\def\crllrl#1{\begin{corollary}\label{#1}}
\def\ecrllr{\end{corollary}}
\def\css{\begin{cases}}
\def\ecss{\end{cases}}

\def\proof{\noindent{\bf Proof}.\ }

\def\fD{{\frak D}}
\def\fE{{\frak E}}
\def\fP{{\frak P}}

\DeclareMathOperator{\aut}{Aut}
\DeclareMathOperator{\alt}{Alt}
\DeclareMathOperator{\AGL}{AGL}

\DeclareMathOperator{\cyc}{cyc}

\DeclareMathOperator{\GF}{GF}
\DeclareMathOperator{\GL}{GL}

\DeclareMathOperator{\hol}{Hol}
\DeclareMathOperator{\id}{id}

\DeclareMathOperator{\Inn}{Inn}

\DeclareMathOperator{\orb}{Orb}

\DeclareMathOperator{\poly}{poly}

\DeclareMathOperator{\PSL}{PSL}
\DeclareMathOperator{\PGL}{PGL}
\DeclareMathOperator{\PGaL}{P\Gamma L}

\DeclareMathOperator{\soc}{Soc}

\DeclareMathOperator{\syl}{Syl}
\DeclareMathOperator{\sym}{Sym}

\newcommand{\fun}{\mathrm{Fun}}

\def\eprf{\hfill$\square$}

\def\qaq{\quad\text{and}\quad}

\newcommand{\grp}[1]{\langle {#1}\rangle}

\begin{document}
\title[Finding a cycle base of a permutation group]{Finding a cycle base 
of a permutation group\\ in polynomial time}
\author{Mikhail Muzychuk}
\address{Netanya Academic College, Netanya, Israel}
\email{muzy@netanya.ac.il}
\author{Ilia Ponomarenko}
\address{Steklov Institute of Mathematics at St. Petersburg, Russia}
\email{inp@pdmi.ras.ru}
\date{}

\begin{abstract}
A cycle base of a permutation group is defined to be a maximal set of its pairwise non-conjugate regular cyclic subgroups. It is proved that a cycle base of a permutation group of degree $n$ can be constructed in polynomial time in~$n$.
\end{abstract}

\date{}
\maketitle

\section{Introduction}
It is well known that the graph isomorphism problem is polynomial-time equivalent to
finding the automorphism group of a graph. However, it is not clear whether the automorphism group given as the input can help to test isomorphism. A byproduct of our main result says that it does help if the input graphs (or any other combinatorial object) is circulant. To be more precise, we need the concept of cyclic base explained below.

Any permutation group $K\le\sym(n)$ acts by conjugation on the set
$$
\cyc(K)=\{G\le K:\ G\ \text{is regular and cyclic}\}.
$$
A {\it cycle base} of $K$ is a set $B\subseteq\cyc(K)$ intersecting each $K$-orbit in this action in exactly one element; in other words, $B$  is a maximal set of pairwise non-conjugate regular cyclic subgroups of~$K$. In slightly different form, this notion was first used in~\cite{P92} for efficient recognizing and isomorphism testing of circulant tournaments.  Then with the help of the classification of finite simple groups, it was proved in~\cite{Mu99} that
$$
|B|\le\varphi(n)
$$
for every cycle base $B$ of the group $K$, where $\varphi$ is the Euler function. Finally, a cycle base technique was applied for polynomial-time recognizing and testing isomorphism of arbitrary circulant graphs \cite{EP03,M04}. In particular, an efficient algorithm was proposed in~\cite{EP03} to construct a cycle base of the automorphism group of a graph.\medskip

The idea to use the cyclic base for Cayley graph isomorphism testing goes back to Babai's lemma~\cite[Lemma~3.1]{B77}, which establishes a one-to-one correspondence between the Cayley representations\footnote{By a Cayley representation of a graph $X$ over a group $G$, we mean a Cayley graph over~$G$ isomorphic to $X$; two such representations are called equivalent if some isomorphism of the corresponding Cayley graphs belong to $\aut(G)$.} of a graph $X$ over group~$G$ and regular subgroups of the group $K=\aut(X)$ that are isomorphic to~$G$. Moreover, two Cayley representations of $X$ are equivalent if and only if the corresponding subgroups are conjugate in~$K$. In this terminology, the above mentioned result~\cite{EP03} shows that if the group $G$ is cyclic, then given the group~$K$, one can efficiently find  a full system of pairwise nonequivalent Cayley representations of the graph~$X$.\medskip

It should be mentioned that not every permutation group $K$ is $2$-closed, i.e., is the automorphism group of a graph; for instance, $\sym(n)$ is a unique 2-transitive group of degree~$n$, which is $2$-closed. In particular,
this may occurs if $K=\aut(X)$, where~$X$ is a combinatorial object defined by a set of relations of arity~$r>2$. Therefore, known algorithms cannot be used to find a full system of pairwise nonequivalent Cayley representations of such $X$ over a cyclic group. The main result of the present paper says that one can find such a system efficiently if the group 
$K$ is given.

\thrml{060616b}
A cycle base of any permutation group of degree $n$ can be constructed
time~$\poly(n)$.
\ethrm

It is assumed that the input permutation group $K$ is given by a set of generators and the cardinality of this set is polynomial in~$n$, see \cite{S03}. The output $B$ is a set of full cycles contained in~$K$; in particular, $B$ is empty if and only if $K$ contains no regular cyclic subgroup. \medskip

For solvable permutation groups, a polynomial-time algorithm for finding cyclic base was constructed in \cite[Theorem~6.3]{EP03}. Therefore, to prove Theorem~\ref{060616b}, it suffices to be able, given a group $K\le\sym(n)$ to find efficiently a solvable group $K_0\le K$ which  {\it controls regular cyclic subgroups} of~$K$, i.e., for each $H\in\cyc(K)$ there exists $k\in K$ such that $H^k\le K_0$ (indeed, in this case, one can find the set $B$ as a subset of cyclic base of~$K_0$, for details, see \cite[Subsection~6.2]{EP03}). The existence of a solvable group $K_0$ which controls regular cyclic subgroups was proven in~\cite{Mu99}, but no algorithm for finding $K_0$ was provided there. In the present paper, we fill this gap in the following theorem, the proof of which occupies the rest of the paper.

\thrml{060616a}
Given a group $K\le\sym(n)$, one can construct in time~$\poly(n)$
a solvable subgroup of $K$ which controls its regular cyclic subgroups.
\ethrm

The Main Algorithm for proving Theorem~\ref{060616a} is described in Section~\ref{260716a}. The basic idea of the algorithm is, as in the case of 2-closed groups, to construct a solvable subgroup by ``removing'' non-abelian composition factors of~$K$ step by step. This reduces the problem to the case, where all non-abelian composition factors are isomorphic to a simple group~$T$ and contained in the socle of~$K$. In contrast to the case of 2-closed groups, where $T$ is an alternating group, here $T$ might be also a projective special linear group. A relevant theory for the "removing" part is developed in Sections~\ref{201016a} and~\ref{041116a}. It is based on the classification of primitive groups containing a regular cyclic subgroup obtained by G.~A.~Jones in~\cite{J02}. This classification is also used in Section~\ref{260716r} providing algorithmic tools to find a regular cyclic subgroup of a primitive groups.\medskip

Let $X$ be a {\it combinatorial object}, i.e., an object in a concrete category in the sense of~\cite{B77}. It is said to be a {\it circulant object} if the group $\aut(X)$ contains a regular cyclic subgroup. A circulant representation of $X$ is defined in the same way as for graphs. Then in view of the one-to-one correspondence between the circulant representations of $X$ and regular cyclic subgroups of~$\aut(X)$, the following statement is an immediate consequence of Theorem~\ref{060616b}.

\crllrl{051116a}
Given a combinatorial object $X$ of size $n$ and the group $\aut(X)$, one can test in time $\poly(n)$ whether $X$ is a circulant object and (if so) find a full system of pairwise nonequivalent circulant representations of $X$ within the same time.
\ecrllr

Using Corollary~\ref{051116a}, it is a routine task to construct a canonical form of a circulant object~$X$ (see e.g. the proof of \cite[Theorem~1.2]{EP03}). 

\crllrl{060616c}
The problem of finding a canonical form of a circulant object~$X$ is polynomial-time reduced to constructing the group $\aut(X)$.
\ecrllr

All undefined notations and standard facts concerning permutation groups can be found in the monographs~\cite{DM} and~\cite{W64}. Throughout the paper, we freely use known polynomial-time algorithms for permutation groups ~\cite[Section~3.1]{S03}.\medskip

{\bf Notation.}

Hereinafter, $\Omega$ denotes a set of cardinality $n$ and $\sym(\Omega)=\sym(n)$ is the symmetric group on $\Omega$.

The orbit set of a group $K\le\sym(\Omega)$ is denoted by $\orb(K)=\orb(K,\Omega)$.

The restriction of the group $K$ to a $K$-invariant set $\Delta\subseteq\Omega$ is denoted by $K^\Delta$.

The pointwise and setwise stabilizers of the set $\Delta$ in the group~$K$  are denoted by $K_\Delta$ and $K_{\{\Delta\}}$, respectively; we also set $K^\Delta=(K_{\{\Delta\}})^\Delta$.

For an imprimitivity system $\fD$ of a group $K\le\sym(\Omega)$, we denote by $K^\fD$ and $K_\fD$ the permutation group induced by the action of $K$ on the blocks of~$\fD$  and the subgroup of $K$ leaving each block of~$\fD$ fixed.


The holomorph $\hol(G)$ of a regular group $G\le\sym(\Omega)$ is identified with the subgroup
of~$\sym(\Omega)$ induced by the right multiplications and automorphisms of~$G$.

The normalizer of $G\leq \sym(\Omega)$ in $\sym(\Omega)$ is denoted by $N_\Omega(G)$.

\section{Finding a cyclic base of a primitive group}\label{260716r}

\sbsnt{Classification.} The primitive groups containing a regular cyclic subgroup were completely described in paper~\cite[Theorem~3]{J02} modulo the classification of finite simple groups. Below, we cite the corresponding result.

\thrml{070616b}
Let $K\le\sym(n)$ be a primitive group containing a regular cyclic subgroup. Then one of the following statements holds:
\nmrt
\tm{1} $C_p\le K\le \AGL_1(p)$, where $n=p$ is prime,
\tm{2} $K=\sym(n)$ for some $n\ge 2$ or $K=\alt(n)$ for some odd $n\ge 3$,
\tm{3} $\PGL_d(q)\le K\le \PGaL_d(q)$ and $n=(q^d-1)/(q-1)$ for some $d\ge 2$,
\tm{4} $K=\PSL_2(11)$, $M_{11}$, or $M_{23}$, and $n=11$, $11$, or $23$,
respectively.
\enmrt
\ethrm

The following auxiliary statement follows from Theorem~\ref{070616b} and will be used in Section~\ref{041116a}. The authors are thankful to Prof. E.~Vdovin for his help with handling the centralizers of graph automorphisms of the projective special linear groups.

\lmml{020716d}
Let $K\le\sym(n)$ be a primitive group and $G\in\cyc(K)$. Then
$$
C_{\aut(K)}(G)\le\Inn(K).
$$ 
\elmm
\proof The group $N:=N_{\sym(n)}(K)$ is embedded into the group $\aut(K)$. Moreover, $N=\aut(K)$ in all cases mentioned in Theorem~\ref{070616b} unless $K=\sym(6)$ or case (3) occurs. In the former case, the statement of lemma is easily checked with the help of~GAP~\cite{gap}. In the remaining case,
$$
\PGL_d(q)\le K\le \PGaL_d(q)
$$ 
and hence the group $\aut(K)$ is embedded into $\aut(\PSL_d(q))$. Then any external automorphism of $K$ can be written as the product of diagonal, field, and graph automorphisms. The automorphisms of the first two types are realized in $\sym(n)$, and one can apply the above argument. Let now $\sigma$ be a graph automorphism of $S$. According to \cite[Table~4.5.1]{GLS},
$$
C_S(\sigma)\in\{S^{}_d(q),O^{}_d(q),O^\pm_d(q)\},
$$
where the parity of~$d$ determines which group on the right-hand side occurs as the centralizer $C_S(\sigma)$. The order of each of these group is not divisible by the Zsigmondy prime for $(q,d)$. Therefore, none of this groups contains an element of order $n=|G|$. This shows that no graph automorphism of $K$ centralizes~$G$.\eprf

\sbsnt{Recognizing.}\label{311016j}
Our first goal is to recognize the groups $K$ appearing in Theorem~\ref{070616b}, and then in each case, to construct a regular cyclic subgroup of~$K$. This is done in more or less standard way in the following statement, where we use the fact that the socle of a subgroup of $\sym(n)$ can be found in polynomial time in~$n$~\cite[Section~3.1]{S03}.

\lmml{080616u}
Given  a primitive group $K\le\sym(n)$, one can test  in time~$\poly(n)$  whether $\cyc(K)\ne\varnothing$, and (if so) find  $H\in\cyc(K)$ within the same time.
\elmm
\proof Let $S=\soc(K)$. If the number $n$ is prime and $S\cong C_n$, then case~(1)  of Theorem~\ref{070616b} occurs and we output $H=S$. Next, if $|S|=(n!)/2$, then $S=\alt(n)$ and case~(2) occurs. Here we output $\cyc(K)=\varnothing$ if $K=S$ and $n$ is even, and the group $H$ generated by a full cycle of $\sym(n)$, otherwise. In the remaining two cases, the group $S$ is determined by its order
up to isomorphism \cite[Theorem~5.1]{KLST}. If now case~(4) occurs, then the group~$H$ can be found by the inspection of the elements of~$K$. Finally, if $S\not\cong\PSL_d(q)$ for some $d$ and $q$, for which $n=(q^d-1)/(q-1)$, then $\cyc(K)=\varnothing$.\medskip

To complete the proof, we assume that $S\cong\PSL_d(q)$ for suitable $d$ and~$q$. Then using the main algorithm from \cite{KS}\footnote{Note that this algorithm is polynomial in $n$, because $q\le n$.}, one can find a $d$-dimensional vector space $V$ over $\GF(q)$ and an explicit isomorphism
$$
f:S\to \PSL(V)
$$ 
given by the images of generators of~$S$. This isomorphism is induced by a bijection from $\{1,\ldots,n\}$ onto the lines of~$V$. Let us extend $f$ to an isomorphism 
$$
f':S'\to\PGL(V),
$$ 
where $S'$ is a unique subgroup of $K$, which is isomorphic to $\PGL_d(q)$ and contains~$S$. Using the natural basis of~$V$, one can construct a Singer subgroup $H\le \GL(V)$ given by an explicit generator matrix. Then $(f')^{-1}(H)$ is a regular cyclic subgroup of~$K$, as required.\eprf\medskip

Let $K\le\sym(\Omega)$ be a primitive group and $H$ a regular cyclic subgroup of~$K$. Denote by $p$   the largest prime divisor of the number $n$. The group $H$ has  a unique subgroup $P=P(H)$ of order~$p$ and the set $\fP=\orb(P,\Omega)$ is an imprimitive system of~$H$. Note that $H$ normalizes the direct sum of the permutation groups $N_\Theta(P^\Theta)$, $\Theta\in\fP$. It is easy to see that the group
\qtnl{311016e}
N(H)=H\,\prod_{\Theta\in\fP} N_\Theta(P^\Theta),
\eqtn 
is permutation isomorphic to the wreath product  $N_\Theta(P^\Theta)\wr H^\fP$ in the imprimitive action. In particular, it is solvable. Note also that the setwise stabilizer of $\Theta$ in $N(H)$ contains~$P$.

\thrml{130616a}
Let $K\le\sym(n)$ be a primitive group containing a regular cyclic subgroup~$H$. Then given $C\in\cyc(K)$, there exists $s\in \soc(K)$ such that 
$$C^s\le N(H).$$
\ethrm
\proof In what follows, we set $P=P(H)$, $N=N(H)$, and $S=\soc(H)$. Assume first that $n=p$ is a prime. Then $P=H$ and $N=N_{\sym(n)}(P)$. Moreover, $\cyc(K)=\syl_p(S)$ and by the Sylow theorem, every $C\in\cyc(K)$ is $S$-conjugate to $P\le N$. This proves the required statement in the considered case. Assume now that $n$ is composite. Then by Theorem~\ref{070616b}, the group $S$ is either projective special or alternating. Let us consider these two cases separately.\medskip

Let $S=\PSL_d(q)$ for appropriate $d$ and~$q$. If $(d,q)=(2,8)$ and $K=\PGaL_2(8)$, then $n=9$, $p=3$, and the required statement follows by a direct calculation in GAP~\cite{gap}. Otherwise, from \cite[Corollary~2]{J02} it follows that
$$
\cyc(K)=\cyc(\PGL_d(q))
$$
and every two groups in $\cyc(K)$ are conjugate in $\PGL_d(q)$. On the other hand, the group $HS\le K$ contains a regular cyclic group $H$ and hence contains the group $\PGL_d(q)$ again by \cite[Corollary~2]{J02}. Thus, every $C\in\cyc(K)$ is conjugate to $H\le N$ in $HS$ and hence in $S$, as required.\medskip

Finally, let $S=\alt(\Omega)$, $h$ a generator of $H$, and $c\in\sym(\Omega)$ a full cycle. Then the cyclic representations of the permutations $h^m$ and $c^m$, where $m=n/p$, consist of $m>1$ cycles of length~$p$. Therefore, these permutations are conjugate in $S$ (see, e.g.~\cite[Lemma~1.2.10]{JK}). Thus, without loss of generality, we may assume that 
$$
\orb(Q,\Omega)=\orb(P,\Omega)=\fP,
$$
where $Q=P(C)$. Thus, the groups $C$ and $H$ can be treated as subgroups of the wreath product $\sym(p)\wr\sym(m)$ in imprimitive action. This product contains the elements 
$$
c'=(1,\ldots,1;\, c_m)\qaq h'=(1,\ldots,1;\, h_m),
$$ 
where $c_m$ and $h_m$ are generators of the groups $C^\fP$ and $G^\fP$, respectively. Each of the permutations $c'$, $h'$ is the disjoint union of $p$ cycles of length~$m$. Therefore, $(c')^{s'}=h'$ for some $s'\in S$. Thus, we may assume that
$$
C\le Q\wr H^\fP\qaq H\le P\wr H^\fP.
$$
If $p=2$, then $Q=P$ and we are done, because in this case $N=P\wr H^\fP$. Let $p>2$. Then $Q$ and $P$ are Sylow subgroups of $\alt(p)$ and hence conjugate in $\alt(p)$. Therefore,  there exists $s\in S$ such that
$$
C^s\le (Q\wr H^\fP)^s=P\wr H^\fP=N,
$$ 
as required.\eprf

\section{Imprimitivity systems of feasible groups}\label{201016a}

\sbsnt{Feasible groups.}\label{271016r}
Let $K\le\sym(\Omega)$ be a transitive group and $\fD$ a minimal imprimitivity system of~$K$. We assume that $K$ is {\it feasible} with respect to $\fD$, by which we mean that  $\fD$ is normal, i.e.,
$$
\orb(K_\fD,\Omega)=\fD,
$$
and $K^\Delta$ is a non-solvable group containing a regular cyclic subgroup for some (and hence for all) $\Delta\in\fD$. In what follows, 
$S=\soc(K_\fD)$.

\lmml{261016u}
In the above notation, $\orb(S,\Omega)=\fD$. Moreover, for each $\Delta\in\fD$, the following statements hold:
\nmrt
\tm{1} $S^\Delta=\soc(K^\Delta)$,
\tm{2} $S^\Delta$ is a $2$-transitive non-abelian simple group.
\enmrt
\elmm
\proof The characteristic subgroup $S$ of the group $K_\fD\trianglelefteq K$ is normal in~$K$. Therefore, $\orb(S,\Omega)$ is a nontrivial imprimitive system of $K$ that is a refinement of $\fD$. Since the imprimitive system is minimal and normal, this implies that  $\orb(S,\Omega)=\orb(K_\fD,\Omega)=\fD$. Next, let $\Delta\in\fD$. Then $\soc(K^\Delta)$ is a simple group (Theorem~\ref{070616b}).  Since $S^\Delta$ is a nontrivial normal subgroup of~$K^\Delta$, we have
$$
\soc(K^\Delta)\le S^\Delta.
$$ 
Since $S^\Delta$ is a direct product of simple groups and $\soc(K^\Delta)$ is normal in $S^\Delta$, we obtain that $S^\Delta = \soc(K^\Delta) C_{S^\Delta}(\soc(K^\Delta))$. It follows from Theorem~\ref{070616b} that $\soc(K^\Delta)$ is a non-abelian $2$-transitive simple group with trivial $C_{\sym(\Delta)}(\soc(K^\Delta)$. Therefore $\soc(K^\Delta) = S^\Delta$ implying both statements (1) and (2).\eprf\medskip

From Theorem~\ref{070616b}, it follows that if a permutation group $K$ contains a regular cyclic subgroup, then $K$ is feasible with respect to every minimal imprimitivity system $\fD$ such that the group $K^\Delta$ is not solvable for some (and hence for all) $\Delta\in\fD$.

\sbsnt{An imprimitivity system induced by stabilizers.}\label{301016d}
Throughout this subsection, $K$ is a feasible group with respect to a minimal imprimitivity system~$\fD$. In the notation of Subsection~\ref{271016r}, we define a binary relation $\sim$ on the set~$\Omega$ by setting 
$$
\alpha\sim\beta\quad\Leftrightarrow\quad S_\alpha\ \text{and}\ S_\beta\ \text{are $S$-conjugate}.
$$
Clearly, this is an equivalence relation. Since $S$ is normal in $K$, the equivalence classes of $\sim$ form an imprimitivity system $\fE$ of $K$.\medskip

From Lemma~\ref{261016u}, it follows that $\fD$ is a refinement of~$\fE$. In a natural way, this induces an equivalence relation on~$\fD$, which is denoted again by~$\sim$; thus for all $\Delta$ and $\Gamma$ in $\fD$, we have $\Delta\sim\Gamma$ if and only if $S_\alpha$ and $S_\beta$ are conjugate in $S$ for all $\alpha,\beta$ belonging to $\Delta\,\cup\,\Gamma$. In the following statement, we identify a bijection with its graph treated as a binary relation.  

\lmml{261016c}
Let $\Delta,\Gamma\in\fD$. Then $\Delta\sim\Gamma$ if and only if there exists a bijection $f:\Delta\to \Gamma$ such that
\qtnl{261016d}
\orb(S,\Delta\times\Gamma)=\{f,f^c\},
\eqtn
where $f^c$ is the complement of $f$ in $\Delta\times\Gamma$. 
\elmm
\proof To implication $\Leftarrow$ in formula~\eqref{261016d} immediately follows from the equality $S_\delta=S_{f(\delta)}$ that holds for all $\delta\in\Delta$ whenever $f\in\orb(S,\Delta\times\Gamma)$. Conversely, assume $\Delta\sim\Gamma$ and take $\delta\in\Delta$. Then there exists
$\gamma\in\Gamma$ such that $S_\delta=S_\gamma$: indeed, since 
the groups $S_{\delta^{}}$ and $S_{\gamma'}$ are $S$-conjugate for all $\gamma'\in\Gamma$, we have
$$
S_{\delta^{}}=(S_{\gamma'})^s=S_{\gamma^{}},
$$
for some $s\in S$, where $\gamma=(\gamma')^s$. Next, the mapping $f:\delta^s\mapsto\gamma^s$, $s\in S$, is well-defined, because if $\delta^s=\delta^t$ for some $t\in S$, then $S_{\gamma^s}=S_{\delta^s}=S_{\delta^t}=S_{\gamma^t}$ and hence $\gamma^s=\gamma^t$ by the $2$-transitivity of $S_\Gamma$ (statement~(2) of Lemma~\ref{261016u}). Similarly, the $2$-transitivity of $S_\Delta$ implies that the mapping $f$ is a bijection. This proves equality~\eqref{261016d}, because from \cite[Corollary~13, p.~86]{W76} it follows that the group $S$ being $2$-transitive on $\Delta$ and $\Gamma$ has at most two orbits in the action on $\Delta\times\Gamma$.\eprf\medskip

For a class $\Lambda\in\fE$, denote by $1_{\fD_\Lambda}$ the identity subgroup on the set 
$\fD_\Lambda$ of all classes $\Gamma\in\fD$ contained in~$\Lambda$.
Given a class $\Delta\in\fD_\Lambda$, we define a bijection 
\qtnl{281016t}
f_\Lambda:\Lambda\to \Delta\times\fD_\Lambda,\quad \gamma\mapsto(f_{\Gamma,\Delta}(\gamma),\Gamma),
\eqtn
where $\Gamma\in\fD$ is a unique block containing $\gamma$ and $f_{\Gamma,\Delta}$ is the inverse to the bijection~$f$ defined in Lemma~\ref{261016c}.

\thrml{281016r}
Let $K$ be a feasible group with respect to a minimal imprimitivity system~$\fD$. Then for each $\Lambda\in\fE$ and each $\Delta\in\fD_\Lambda$, the bijection~\eqref{281016t} induces a permutation isomorphism from $(K_\fD)^\Lambda$ onto the direct product $(K_\fD)^\Delta\times 1_{\fD_{\Lambda}}$. 
\ethrm
\proof Note that the set  $\Delta$ is of cardinality at least five by statement~(2) of Lemma~\ref{261016u}. In the notation of Lemma~\ref{261016c}, this immediately implies that $|f|<|f^c|$. Therefore,  the bijection $f$ is $K_\fD$-invariant by the normality of the group~$S$ in $K_\fD$. It follows that for each $\Gamma\in\fD_\Lambda$, all $k\in K_\fD$ and $\gamma\in\Gamma$, we have
$$
(\gamma^k)^{f_\Lambda}=
(f(\gamma^k),\Gamma)=
(f(\gamma)^k,\Gamma)=
(f(\gamma),\Gamma)^{(k,1)}=
(\gamma^{f_\Lambda})^{(k,1)},
$$
where $f=f_{\Gamma,\Delta}$. This proves that $f$ is a permutation isomorphism from $(K_\fD)^\Lambda$ onto $(K_\fD)^\Delta\times 1_{\fD_\Lambda}$.\eprf

\sbsnt{An imprimitivity system associated with $S$.}\label{301016e}
The group $S$ is the direct product of, say $d$, simple groups $S_i$, $i=1,\ldots,d$. By Lemma~\ref{261016u}, every set $\Delta\in\fD$ is $S_i$-invariant. Therefore the group $(S_i)^\Delta$ is either trivial or isomorphic to $S_i$. Set $\fE'=\{\Omega_1,\ldots,\Omega_d\}$, where
$$
\Omega_i=\bigcup_{(S_i)^\Delta\ne 1}\Delta.
$$
It is easy to see that $\fD$ is a refinement of $\fE'$ and the restriction of $S_i$ to $\Omega_i$ is isomorphic to~$S_i$. Furthermore, by statement~(2) of Lemma~\ref{261016u} the sets $\Omega_i$ are pairwise disjoint. Consequently, $\fE'$ is a partition of~$\Omega$. This shows that
\qtnl{291016a}
S=\prod_{i=1}^d(S_i)^{\Omega_i}.
\eqtn
Note that by the transitivity of $K$ and statement~(1) of Lemma~\ref{261016u}, the groups $\soc(S^\Delta)$ with $\Delta\in\fD$ are pairwise isomorphic. Since $S_i\cong (S_i)^\Delta$, this implies that the groups $S_i$ are pairwise isomorphic too.

\lmml{290616e}
Let $K$ be a feasible group with respect to the imrimitivity system~$\fD$. Then the set $\fE'$ is an imprimitivity system of~$K$. Moreover, $\fE$ is a refinement of~$\fE'$. 
\elmm  
\proof The simple groups $S_1,\ldots,S_d$ are uniquely determined up to permutation of indices. Therefore, the group~$K$ acts on the set of all of them by conjugation and hence permutes the sets~$\Omega_i$. This proves the first statement of the lemma. Next, in view of formula~\eqref{291016a} for each $i$ and $\delta\in\Omega_i$, we have
\qtnl{291016v}
S_\delta=(S_i)_{\delta}\prod_{j\ne i}S_j.
\eqtn
Since the $S_i$ are isomorphic non-abelian simple groups, this implies 
that the groups $S_\delta$ and $S_\gamma$ are not $S$-conjugate unless $\gamma\in\Omega_i$. Therefore,  the class of $\fE$ containing $\delta$ is a subset of $\Omega_i$. Thus, $\fE$ is a refinement of $\fE'$. \eprf

\rmrkl{301016a}
From Lemma~\ref{290616e} and formula~\ref{290616e}, it immediately follows that $\fE=\fE'$ if and only if $(S_i)_\delta=(S_i)_\gamma$ for all $i=1,\ldots,d$ and all $\delta,\gamma\in \Omega_i$.
\ermrk

\section{Imprimitive groups containing regular cyclic subgroups}\label{041116a} 

\sbsnt{The imprimitivity systems $\fE$ and $\fE'$ are equal.} In this subsection, we apply the theory developed in Section~\ref{201016a} to feasible groups containing a regular cyclic subgroup. 

\thrml{301016b}
Let $K\le\sym(\Omega)$ be a feasible group with respect to a minimal imprimitivity system~$\fD$. Suppose that $K$ contains a regular cyclic subgroup. Then the imprimitivity systems $\fE$ and $\fE'$ defined in Subsections~\ref{301016d} and~\ref{301016e} coincide.
\ethrm

Let $C$ be a regular cyclic subgroup of the group~$K$. Then in the notation of Section~\ref{201016a}, the group $(S_i)^{\Omega_i}C^{\Omega_i}$ satisfies the hypothesis of Theorem~\ref{201016e} below. Therefore in view of Remark~\ref{301016a}, the statement of Theorem~\ref{301016b} is an immediately consequence of Theorem~\ref{201016e}, the proof of which occupies the rest of this subsection.

\thrml{201016e}
Let $K\le\sym(\Omega)$ be a group containing a regular cyclic subgroup and  $S\trianglelefteq K$ a non-abelian simple group. Assume that
the imprimitivity system
$$
\fD=\orb(S,\Omega)
$$
of the group $K$ is minimal. Then for any  $\alpha,\beta\in\Omega$, the groups $S_\alpha$ and $S_\beta$ are conjugate in~$S$. 
\ethrm
\proof  To prove the first statement, without loss of generality we may assume that $K=SC$, where $C$ is a regular cyclic subgroup of~$K$. It is  easy to see that in this case  
\qtnl{201016g}
K_\fD=SC_\fD.
\eqtn
and the cyclic group $C^\fD$  is regular. Note that by statement~(2) of Lemma~\ref{261016u}, for any set $\Delta\in\fD$, the group $(K_\fD)^\Delta\ge S^\Delta$ is primitive and contains a regular subgroup $C^\Delta$.\medskip

Recall that the group $S$ being simple acts on a set $\Delta\in\fD$ faithfully. Therefore, the intersection of $S$ with the pointwise stabilizer $K_\Delta$ of the set $\Delta$ in the group~$K$ is trivial. Since $K=SC$, this shows that $K_\Delta\trianglelefteq K_\fD$ is a cyclic group. Consequently, 
$$
H=\grp{(K_\Delta)^c:\ c\in C}
$$ 
is a normal solvable subgroup of~$K_\fD$. Therefore, $H$ intersects $S$ trivially and hence is cyclic. It follows that each group $(K_\Delta)^c\le H$ is also cyclic. Since all these groups are of the same order, they must be equal. Thus, the group $K_\Delta$ and  hence the group $(K_\fD)_\Delta$ fixes each point of the set $\Omega$. This proves the following statement.

\lmml{201016t}
For each $\Delta\in\fD$, the restriction epimorphism $\pi:K_\fD\to (K_\fD)^\Delta$ is an isomorphism.
\elmm

Each $c\in C$ induces the automorphism $\sigma_c:k\mapsto k^c$ of the group $K_\fD$ that centralizes the group $C_\fD$.  Therefore  $\pi^{-1}\sigma_c\pi$ is an automorphism of the primitive group $(K_\fD)^\Delta$ that centralizes the regular cyclic subgroup $(C_\fD)^\Delta$ (Lemma~\ref{201016t}). By Lemma~\ref{020716d}, this implies that $\pi^{-1}\sigma_c\pi$ is an inner automorphism of the group $(K_\fD)^\Delta$ corresponding to a certain element $k'\in (K_\fD)^\Delta$. Thus, $\sigma_c$ equals the inner automorphism of $K_\fD$ corresponding to the element $k=\pi^{-1}(k')$.

\lmml{201016u}
For each $c\in C$, there exists $k\in K_\fD$ such that $x^c=x^k$ for all $x\in K_\fD$.
\elmm

To complete the proof of the first statement of Theorem~\ref{201016e}, let $\alpha,\beta\in\Omega$. By the transitivity of $C$, there exists $c\in C$ such that $\beta=\alpha^c$. In view of Lemma~\ref{201016u}, one can find $k\in K_\fD$, for which
\qtnl{300616r}
S_\beta=S_{\alpha^c}=(S_\alpha)^c=(S_\alpha)^k=S_{\alpha^k}.
\eqtn
Note that by the definition of the group~$K_\fD$, the points $\alpha$ and $\alpha^k$ belongs to the same orbit of the group~$S$. Therefore, $\alpha^k=\alpha^s$ for some $s\in S$. Thus,  
$$
S_{\alpha^k}=S_{\alpha^s}=(S_\alpha)^s,
$$
which together with~\eqref{300616r} shows that $S_\beta=(S_\alpha)^s$, as required.\eprf

\sbsnt{The embedding into the wreath product.}
Under the hypothesis of Theorem~\ref{301016b}, we fix a $K$-block $\Delta\in\fD$ and arbitrary elements $k_\Lambda\in K$ taking the class $\Lambda_\Delta\in\fE$ containing $\Delta$ to the class $\Lambda\in\fE$ (here we assume that $k_{\Lambda_\Delta}=1$). Let us define a bijection
\qtnl{301016z}
f^*:\Omega\to \Delta\times\fD,\ \gamma\mapsto(\gamma^*,\Gamma)
\eqtn
with $\gamma^*=\gamma^{k_{\Lambda_{}} f_{\Lambda_\Delta}}$, where $\Lambda$ is the class of $\fE$ that contains $\gamma$ (and hence $\Gamma$) and $f_{\Lambda_\Delta}$ is the bijection defined in formula~\eqref{281016t} for $\Lambda=\Lambda_\Delta$. Thus if $K^*=K^{f^*}$, then
\qtnl{311016y}
K^*\le \sym(\Delta)\wr\sym(\fD),
\eqtn
where the wreath product on the right-hand side is considered in the imprimitive action.

\thrml{221016b}
Under the identification of $K$ and $K^*$ via the bijection~$f^*$, the following statements hold:
\nmrt
\tm{1} $K_\fD$ is the direct sum ot the permutation groups $(K_\fD)^\Lambda$,
$\Lambda\in\fE$,
\tm{2} $(K_\fD)^\Lambda=(K_\fD)^\Delta\times 1_{\fD_\Lambda}$ for all $\Lambda\in\fE$.
\enmrt
\ethrm
\proof Statement~(1) follows from Theorem~\ref{301016b} and formula~\eqref{291016a}, whereas statement~(2) is a straightforward consequence of Theorem~\ref{281016r}.\eprf\medskip

Let $H$ be a regular cyclic subgroup of the group $K^\Delta$. 
The wreath product on the right-hand side of inclusion~\eqref{311016y} contains the subgroup
\qtnl{301016u}
W^*=W^*(H,\Delta)=N^*(H)\wr K^\fD
\eqtn
where the group $N^*(H)$ is defined  as follows. Let the group $P=P(H)$ and the imprimitivity system $\fP$ be as in Subsection~\ref{311016j}.  Then 
\qtnl{021116a}
N^*(H)=N_\Omega(H)\,\prod_{\Theta\in\fP} N_\Theta(P^\Theta).
\eqtn
Note that this group contains the group $N(H)$ defined in~\eqref{311016e}.
The following statement is crucial for our arguments.

\thrml{160616t}
Let $K\le\sym(\Omega)$ be a transitive group and $\Delta$ a minimal $K$-block. Suppose that $K^\Delta$ is a non-solvable group containing a regular cyclic subgroup $H$. Then 
\nmrt
\tm{1} $K$ is a feasible group with respect to the imprimitivity system $\fD=\Delta^K$,
\tm{2} the group $K^*\cap W^*$ controls the regular cyclic subgroups of $K^*$.
\enmrt
\ethrm
\proof The first statement follows from the definition. To prove the second one, we identify $K$ and $K^*$ via the bijection~$f^*$. In what follows, we assume that the imprimitivity system $\fE$ consists of $d\ge 1$ blocks, say $\Lambda_1,\ldots,\Lambda_d$. The number $e$ of the  blocks of $\fD$ contained in $\Lambda_i$ does not depend on $i$; these blocks are denoted by $\Delta_{i1},\ldots,\Delta_{ie}$. In this notations, $|\fD|=de$ and $\Delta_{ij}=\Delta$.\medskip

It suffices to verify that for every regular cyclic subgroup $C\le K$ there exists an element $k\in K$ such that $C^k\le W^*$. To this end, we make use of Theorem~\ref{221016b} to write a generator $c$ of the group $C_\fD\le K_\fD$ in the form
\qtnl{170616a}
c=(c^{\Lambda_1},\ldots,c^{\Lambda_d})=
(\underbrace{c_1,\ldots,c_1}_e,\cdots,\underbrace{c_d,\ldots,c_d}_e)
\eqtn
where $c_i\in (K_\fD)^\Delta$ for $i=1,\ldots,d$. Note that $K^\Delta$ is a primitive group containing a cyclic regular subgroup~$H$. Therefore by Theorem~\ref{130616a}, one can find an element $s_i\in S^\Delta$ such that 
\qtnl{200616a}
(c_i)^{s_i}\le N(H),\qquad i=1,\ldots,d,
\eqtn
where $N(H)$ is the group defined by formula~\eqref{311016e}.
By statement~(1) of Theorem~\ref{221016b}, the permutation $s\in\sym(\Omega)$ such that $s^{\Delta_{ij}}=s_i$ for all $i,j$, belongs to the group $S\le K$. In particular, $c^s\in K$. Together with formula~\eqref{200616a}, this shows that
\qtnl{021116b}
(C^s)_\fD\le\underbrace{N(H)\times\cdots\times N(H)}_{de}\qaq (C^s)^\fD\le K^\fD.
\eqtn
At this point, we make use of an obvious permutation isomorphism from the group $N^*(H)$ onto the wreath product $N_\Theta(P^\Theta)\wr N_{\fP}(H^\fP)$ with a fixed $\Theta\in\fP$ to identify the set $\Delta$ with $\Theta\times\fP$. Then by the associativity of the wreath product, we have
\qtnl{021116t}
W^*=N_\Theta(P^\Theta)\,\wr\, (N_{\fP}(H^\fP)\wr K^\fD).
\eqtn
Note that by the first inclusion in~\eqref{021116b}, the set $\Theta$ is a block of the group $C^s$; denote by $\fP^*$ the corresponding imprimitivity system.\medskip 

To complete the proof of Theorem~\ref{160616t}, we use Lemma~\ref{011116a} and Corollary~\ref{250117t} proved in the Appendix (Section~\ref{250117y}).
Namely, in view of the second inclusion of~\eqref{021116b}, the groups 
$$
A'=N_\fP(H^\fP),\quad B'=K^\fD,\quad C'=(C^s)^{\fP^*}
$$
and the sets $\Delta'=\fP$ and $\Gamma'=\fD$ satisfy the hypothesis of Lemma~\ref{011116a}. Next, by the definition of $N(H)$, the group $C_0=(C^s)^\fP=H^\fP$ is regular and cyclic. Therefore, 
$$
N_{\fP}(C_0)=N_{\fP}(H^\fP)=A'.
$$
Thus the condition of Corollary~\ref{250117t} follows from the first inclusion in~\eqref{021116b}. Thus, by this corollary we obtain
\qtnl{021116y}
(C^s)^{\fP^*}\le N_{\fP}(H^\fP)\wr K^\fD.
\eqtn
Now we again apply Corollary~\ref{250117t} but this time to the groups
$$
A'=N_\Theta(P^\Theta),\quad B'=N_\Theta(H^\fP)\wr K^\fD,\quad C'=C^s
$$
and the sets $\Delta'=\Theta$ and $\Gamma'=\fP^*$. Note that the hypothesis of Lemma~\ref{011116a} follows from the definition of $\fP^*$ and formula~\eqref{021116y}, respectively. The first inclusion ~\eqref{021116b} and the definition of~$N(H)$ imply that $(C^s)^\Theta$ is a regular cyclic subgroup of the group 
$$
N(H)^\Theta=N_\Theta (P^\Theta)\cong\AGL(1,p)
$$
Since $p=|\Theta|$ is a prime, this subgroup is unique and hence the condition of Corollary~\ref{250117t} is also satisfied. Thus, by this corollary and formula~\eqref{021116t} we conclude that
$$
G^s\le A'\wr B'=N_{\Theta}(P^\Theta)\wr (N_{\fP}(H^\fP)\wr K^\fD)=W^*
$$
as required.\eprf

\section{The main algorithm and poof of Theorem~\ref{060616a}}\label{260716a}
The Main Algorithm below finds a solvable subgroup $M$ of a given permutation group~$K$ that controls its regular cyclic subgroups. In particular, $\cyc(M)=\varnothing$, whenever $\cyc(K)=\varnothing$. Except for the algorithm constructed in the proof of Lemma~\ref{080616u}, we use the standard algorithms for computing with permutation groups \cite[Section~3.1]{S03} and  the algorithm in \cite[Corollary~6.4]{L93} finding the intersection of two groups in $\sym(n)$ in time~$\poly(n)$, whenever one of them is solvable.\medskip

\centerline{\bf Main Algorithm.}\medskip

\noindent {\bf Input:} a transitive permutation group $K\le\sym(n)$.

\noindent {\bf Output:} a solvable group $M\le K$ that controls the regular cyclic subgroups~$K$.

\medskip

\noindent{\bf Step 1.} If $n=1$, then output $M=K$. Find a minimal $K$-block~$\Delta$ and the imprimitive system  $\fD$ containing $\Delta$. If $\fD$ is not normal, then output $M=\{\id_\Omega\}$.\medskip
	
\noindent{\bf Step 2.}  Recursively apply the algorithm to the group $K^{\fD}\le\sym(\fD)$; replace~$K$ by the full preimage in $\sym(\Delta)$ of the resulting group.\medskip

\noindent{\bf Step 3.} 
If the group $K$ is solvable or intransitive,  then output $M=K$ or $\{\id_\Omega\}$, respectively.\medskip

\noindent{\bf Step 4.} Apply Lemma~\ref{080616u} to check whether  $K^\Delta$ contains a regular cyclic subgroup~$H$. If there is no such $H$, then output  $M=\{\id_\Omega\}$.\medskip

\noindent{\bf Step 5.} Find the bijection $f^*:\Omega\to\Delta\times\fD$ and the group $W^*=W^*(\Delta,H)$, defined by formulas~\eqref{301016z} and~\eqref{301016u}, respectively.\medskip

\noindent{\bf Step 6.} Output the full $f^*$-preimage $M$ of the group $K^*\cap W^*$, where $K^*=K^{f^*}$.\medskip

Let us prove the correctness of the algorithm. The output at Step~1 is obviously correct. After Step~2, we may assume by induction on~$n$ that $K^\fD$ is a solvable group such that if $G$ is a regular cyclic subgroup of the input group, then 
$$
(G^\fD)^{k^\fD}\in\cyc(K^\fD)
$$
for some element $k$ of the input group. Therefore, $G^k\in\cyc(K)$. This means that the group $K$ controls the regular cyclic subgroups of the input group. Thus, the output at Step~3 is correct. The correctness of Step~4 follows from the fact that 
$$
G^\Delta\in\cyc(K^\Delta)\quad\text{for all}\ G\in\cyc(K)\ \,\text{and}\ \, \Delta\in\fD.
$$

\lmml{311016a}
At Step~5, the imprimitivity system $\fD$ of the group $K$ is still minimal and $K$ is a feasible group with respect to~$\fD$. 
\elmm
\proof Indeed, the group $\cyc(K^\Delta)\ne\varnothing$ by Step~4. Moreover, $K^\Delta$ is not solvable, for otherwise $K$ is isomorphic to a subgroup of a solvable group $K^\Delta\wr K^\fD$ and hence is solvable in contrast to Step~3. This implies that the group 
$$
(K_0)^\Delta\ge K^\Delta
$$ 
is also non-solvable, where $K_0$ denotes the input group. Therefore,  $K_0$ is a feasible group with respect to~$\fD$. By Lemma~\ref{261016u}, this implies that 
$$
\fD=\orb(S_0,\Omega), 
$$
where $S_0=\soc((K_0)_{\fD})$. Since obviously $S_0\le K$, the group $K^\Delta\ge (S_0)^\Delta$ is primitive by statement~(2) of that lemma. Thus $\fD$ is a minimal imprimitive system of~$K$.\eprf\medskip

By Lemma~\ref{311016a} and statement~(2) of Theorem~\ref{160616t}, the group $K^*\cap W^*$ constructed at Step~6 controls the regular cyclic subgroups of~$K^*$. Moreover, it is solvable as a subgroup of a solvable group $W^*=N^*(H)\wr K^\fD$: indeed, $N^*(H)$ is solvable by its construction, whereas $K^\fD$ is solvable by the induction. Thus, the ouput group $M\le K$ at Step~6  is solvable and controls the regular cyclic subgroups~$K$. The correctness of the Main Algorithm is completely proved.\medskip

To estimate the running time $f(n)$ of the Main Agorithm, we note that
Steps~1 and~3 run in polynomial time in~$n$. The same is true for Step~4 by  Lemma~\ref{080616u}.  Finding the bijection $f^*$ at Step~5 is reduced to finding the orbits of $S=\soc(K_\fD)$ on the set $\Omega\times\Omega$ (Lemma~\ref{261016c}) and hence can efficiently be implemented in time $\poly(n)$. Finally, within the same time one can find the intersection of $K^*$ and solvable group $W^*$ (see the remark before the Main Algorithm).
Since Step~2 can easily be implemented in time $f(n/m)+n^c$, where
$m=|\Delta|$ is a divisor of $n$ and $c$ is a constant, we get
$$
f(n)\le f(n/m)+n^c.
$$
Taking into account that $m\ge 2$, we conclude that $f(n)=n^{O(1)}$, as required.\eprf

\section{Appendix. Lemma on the wreath product}\label{250117y}
In this section, we establish a sufficient condition for a permutation group to be a subgroup of a wreath product in the imprimitive action. Below given sets $\Delta$ and~$\Gamma$, we denote by $\fD(\Delta,\Gamma)$ the partition of $\Delta\times\Gamma$ into the subsets 
\qtnl{250117a}
\Delta_\gamma=\Delta\times\{\gamma\},\qquad\gamma\in\Gamma,
\eqtn
which are identified with $\Delta$ with the help of the bijection $(\delta,\gamma)\mapsto\delta$. A permutation $g\in\sym(\Delta\times\Gamma)$ belongs to the wreath product $\sym(\Delta)\wr\sym(\Gamma)$ if and only if it preserves the partition $\fD(\Delta,\Gamma)$. In this case, $g$ permutes the blocks~\eqref{250117a} via 
$$
(\Delta_\gamma)^g = \Delta_{\gamma^g}.
$$ 
For each $\gamma\in\Gamma$, the permutation $g$ induces a permutation $g(\gamma)\in\sym(\Delta)$ such that $(\delta,\gamma)^g=(\delta^{g(\gamma)},\gamma^g)$. From the definition, it immediately follows  that 
\qtnl{eq:wr_product}
(f^{-1})(\gamma^f)=f(\gamma)^{-1}\qaq	(fg)(\gamma)=f(\gamma)g(\gamma^{f})
\eqtn
for all $f,g\in\sym(\Delta)\wr\sym(\Gamma)$. In what follows, we set $C(\gamma)=\{c(\gamma):\ c\in C\}$ for all $C\subseteq\sym(\Delta)\wr\sym(\Gamma)$ and $\gamma\in\Gamma$.\medskip

The well-known Kaloujnine-Krasner embedding Theorem \cite[Theorem~2.6A]{DM} implies that an imprimitive permutation group is embedded into a wreath product of its block restriction and the quotient. The statement below gives necessary and sufficient conditions for a permutation group to be contained in a given wreath product of two permutation groups. In the sequel to avoid a confusion, we denote by $\fun(X,Y)$ (not by $Y^X$) the set of all functions from $X$ to $Y$.

\lmml{011116a}
Let $A\le\sym(\Delta)$ and $B\le\sym(\Gamma)$, and let $C$ be a transitive subgroup of $ \sym(\Delta\times\Gamma)$. Suppose that $\fD=\fD(\Delta,\Gamma)$ is an imprimitivity system of~$C$ and $C^\fD\le B$. Then
\nmrt
\tm{1} if $C\le A\wr B$, then  for each $\gamma\in\Gamma$, there exist $t\in\fun(\Gamma,A)$ and $C_0\leq A$ such that $(tCt^{-1})^{\Delta_\gamma} = C_0\times 1_{\{\gamma\}}$;
\tm{2} if $C_0\leq A$ is such that $C^{\Delta_\gamma} = C_0\times 1_{\{\gamma\}}$ for all $\gamma$, then 
$t Ct^{-1}\leq C_0\wr B$ for some $t\in\fun(\Gamma,N_\Delta(C_0))$. 
\enmrt
\elmm
\proof 
Throughout the proof, we fix an arbitrary $\gamma_0\in \Gamma$. By the transitivity of the action of $C$ on $\Gamma$, for each $\gamma\in\Gamma$ one can find  $c_\gamma\in C$ such that $\gamma_0^{c_\gamma} = \gamma$. We assume that $c_{\gamma_0}=1_{\Delta\times\Gamma}$. Set $t$ to be a permutation in $\fun(\Gamma,\sym(\Delta))$ such that 
\qtnl{250117j}
t(\gamma):=c_\gamma(\gamma_0)\quad\text{for all}\ \,\gamma\in\Gamma.
\eqtn

To prove statement~(1), assume that $C\le A\wr B$. Then $c(\gamma)\in A$ for all $c\in C$ and $\gamma\in\Gamma$. Therefore the permutation $t$ defined in~\eqref{250117j} belongs to $\fun(\Gamma,A)$. 
To show that $t$ is the required element, we first note that $C_\gamma(\gamma)$ is a subgroup of the group~$A$ and $C^{\Delta_\gamma} = C_\gamma(\gamma)\times 1_{\{\gamma\}}$, where $C_\gamma=C_{\{\Delta_\gamma\}}$. Now, from the obvious equality $c_\gamma^{-1} C_0 c_\gamma = C_\gamma$, where $C_0=C_{\gamma_0}$. Therefore, by the second formula in~\eqref{eq:wr_product}, we obtain
$$
C_0(\gamma_0)c_\gamma(\gamma_0)=(C_0c_\gamma)(\gamma_0)=(c_\gamma C_\gamma)
(\gamma_0)=c_\gamma(\gamma_0)C_\gamma(\gamma_0^{c_\gamma})
$$
and hence
$$ 
C_0(\gamma_0)^{c_\gamma(\gamma_0)}=C_\gamma(\gamma).
$$
Thus,
$$
((t C t^{-1})_\gamma)(\gamma) =
t(\gamma) C_\gamma(\gamma) t^{-1}(\gamma) =  c_\gamma(\gamma_0) C_\gamma(\gamma) c_\gamma(\gamma_0)^{-1} 
= C_0(\gamma_0) = C_0,
$$
as required.\medskip

To prove statement~(2), assume that $C^\fD\le B$. Since $\fD$ is an imprimitive system of~$C$, this implies that $C\le \sym(\Delta)\wr\sym(\Gamma)$ and the permutations $c(\gamma)\in\sym(\Delta)$ are well-defined for every $c\in C$. For each $d\in C$ and $\gamma\in\Gamma$, formulas~\eqref{eq:wr_product} imply that
\qtnl{eq:tdt}
\begin{array}{c}
	(tdt^{-1})(\gamma)=t(\gamma)d(\gamma^t)t^{-1}(\gamma^{td}) = c_\gamma(\gamma_0) d(\gamma) t^{-1}(\gamma^d) = \\
	c_\gamma(\gamma_0) d(\gamma) t((\gamma^d)^{-1})=
	c_\gamma(\gamma_0) d(\gamma) (c_{\gamma^d}(\gamma_0))^{-1} = 
	(c_\gamma d (c_{\gamma^d})^{-1})(\gamma_0),
\end{array}
\eqtn
where $t$ is the permutation defined in~\eqref{250117j}.
However, $\gamma_0^{c_\gamma d (c_{\gamma^d})^{-1}} = \gamma_0$. So by the assumption of statement~(2),
$$
(c_\gamma d (c_{\gamma^d})^{-1})(\gamma_0)\in C_{\{\Delta_{\gamma_0}\}}(\gamma_0)=C_0.
$$
This together with~\eqref{eq:tdt} implies that $(tdt^{-1})(\gamma)\in C_0$.
Since $(tdt^{-1})\in C^\fD\le B$, we conclude that $tdt^{-1}\in C_0\wr B$. It remains to show that $t\in\fun(\Gamma,N_\Delta(C_0))$, or, equivalently, that $c_\gamma(\gamma_0)\in N_\Delta(C_0)$ for each $\gamma\in\Gamma$. To this end, take arbitrary $\gamma$. Then by the assumption of statement~(2),
given $f\in C_0$, 
$$
f\times 1_{\{\gamma\}}\in C^{\Delta_\gamma}
$$
and hence there exists $d\in C_{\{\Delta_\gamma\}}$ such that $d(\gamma)=f$. Note that $\gamma^d=\gamma$, and also $(tdt^{-1})(\gamma)\in C_0$. Thus, by formula~\eqref{eq:tdt}, we conclude that 
$$
c_\gamma(\gamma_0)\,f\,c_{\gamma}(\gamma_0)^{-1}=c_\gamma(\gamma_0)\, d(\gamma)\,c_{\gamma^d}(\gamma_0)^{-1}\in C_0,
$$
whence $c_\gamma(\gamma_0)\in N_\Delta(C_0)$, as desired.\eprf

\crllrl{250117t}
Under the hypothesis of statement~(2) of Lemma~\ref{011116a}, assume 
$N_\Delta(C_0)\leq A$. Then $C\le A\wr B$.
\ecrllr

\end{document}